\documentclass[11pt]{article}
\usepackage{amsfonts}
\usepackage{amssymb}
\usepackage{graphicx}
\usepackage{amsmath}
\def\Z{{\mathbb Z}}

\def\R{{\mathbb R}}

\def\SS{{\mathbb S}}

\def\bb{\begin}
\def\bc{\begin{center}}       \def\ec{\end{center}}
\def\be{\begin{equation}}     \def\ee{\end{equation}}
\def\ba{\begin{array}}        \def\ea{\end{array}}
\def\bea{\begin{eqnarray}}    \def\eea{\end{eqnarray}}
\def\beaa{\begin{eqnarray*}}  \def\eeaa{\end{eqnarray*}}
\def\hh{\!\!\!\!}             
\def\EQ{\hh & = & \hh}

\def\a{\alpha}                
\def\e{\varepsilon}

               \def\Ga{\Gamma}
            
\def\oo{\infty}               \def\d{\cdot}
               
\def\lb{\label}               
\def\q{\quad}                 \def\qq{\qquad}
\def\f{\frac}                 
\def\z{\left}                 \def\y{\right}
\def\tm{\times}               

\def\A{{\mathcal A}}
\def\B{{\mathcal B}}

\def\F{{\mathcal F}}
\def\P{{\mathcal P}}
\def\S{{\mathcal S}}
\def\L{{\mathcal L}}
\def\M{{\mathcal M}}

\def\apf{\A\P\F}
\def\aps{\A\P\S}
\def\apl{\A\P\L}

\def\vr{\varrho}

\def\stp{\SS_{2\pi}}
\def\dmu{\,{\rm d}\mu}
\def\dnu{\,{\rm d}\nu}
\def\q{\quad}
\def\qq{\qquad}

\def\Minv{\M_{\rm inv}}

\def\Proof{\noindent{\bf Proof} \quad}
\def\qed{\hfill $\Box$ \smallskip}

\def\bu{$\bullet$\ }
\def\z{\left}
\def\y{\right}
\def\disp{\displaystyle}

\def\IR{\in \R}

\def\lb{\label}

\def\ifl{\iffalse}
\def\Proof{\noindent{\bf Proof} \quad}
\def\qed{\hfill $\Box$ \smallskip}
\def\wx#1#2#3#4#5#6#7#8{\bibitem{#1} {#2}, {#3}, {\it #4}, {{\bf #5}} {(#6)},
{#7}--{#8}.}

\begin{document}

\title{\bf Ergodic behaviour of nonconventional ergodic averages for commuting transformations}

\author{\q Xia Pan,$^1$ \q Zuohuan Zheng$^1$ \q and \q Zhe Zhou$^1$ \footnote{All supported by the National Natural Science Foundation of China(NSFC) under grant 11671382}}

\date{}%

\maketitle

\begin{center}

$^1$ Academy of Mathematics and System
Sciences, Chinese Academy of Sciences, Beijing 100190, China

E-mail: {\tt panxia13@mails.ucas.ac.cn} (Xia Pan) \\
{\tt zhzheng@amt.ac.cn} (Z. Zheng)\\
{\tt zzhou@amss.ac.cn} (Z. Zhou)
\end{center}

\bb{center} \today \end{center}
\begin{abstract}
Based on T.Tao's result of norm convergence of multiple ergodic averages for commuting transformation, we obtain there is a subsequence which converges almost everywhere. Meanwhile, the ergodic behaviour, which the time average is equal to the space average, of diagonal measures is obtained and we give different result according to the classification of transformations. Additionally, on the torus $\mathbb{T}^{d}$ with special rotation, such as $R_{\alpha_{1},\cdots,\alpha_{d}}:\mathbb{T}^{d}\rightarrow \mathbb{T}^{d},$ with $1,\alpha_{1},\cdots,\alpha_{d}$  are rationally independent. we can not only get the convergence in T.Tao's paper for every point in $\mathbb{T}^{d}$, but also get a beautiful result for ergodic behaviour.
\end{abstract}

{\small


}


\section{Introduction}
\setcounter{equation}{0}

In 2008, T.Tao proved a convergence result for several commuting transformations :

\bb{thm} {\rm \cite{Tao08}} \lb{tao} Let $d\geq1$ be an integer. Assuming that $T_{1},T_{2},\dots,T_{d}:X\to X$ are
commuting invertible measure-preserving transformations of a measure space $(X,\mathcal{B},\mu)$, Then, for any $f_{1},f_{2},\dots,f_{d}\in L^{\infty}(X,\mathcal{B},\mu)$, the averages
\be \lb{1.1} \lim_{N \to \oo}\disp \f 1 N \sum_{n=0}^{N-1} f_{1}(T_{1}^{n}x)\cdots f_{d}(T_{d}^{n}x)\ee
are convergent in $L^{2}(X,\mathcal{B},\mu)$.
\end{thm}

Soon after, H.towsner\cite{Tow09}, B.Host\cite{BH09} and T.Austin\cite{Aus10} gave different proofs of Theorem \ref{tao}. T.Tao's approach was combinatorial and finitary, inspired by the hypergraph regularity and removal lemmas. H.towsner used nonstandard analysis. T.Austin and B.Host all gave proofs of the same result by ergodic methods, building an extension of the original system with good properties.

There are a rich history towards Theorem \ref{tao}. For $d=1$, it reduces to the classical mean ergodic theorem. When $T_{1}=T,T_{2}=T^{2},\cdots,T_{d}=T^{d}$, Furstenberg studied such averages originally in his proof of Szemer\'{e}di's theorem\cite{F77}, where $T$ is weakly mixing or for general $T$ but $d=2$. For higher $d$, various special cases have been shown by Conze and Lesigne\cite{CL84},\cite{CL87}, Furstenberg and Weiss\cite{FW96}, Host and Kra\cite{HK01}, and Ziegler\cite{Zie02}. Finally, it was totally proved by Host and Kra\cite{HK05} for arbitrary $d$, and independently by Ziegler\cite{Zie07}.

When $T_{1},T_{2},\dots,T_{d}$ are commuting measure-preserving transformations. With some hypothesis on the transformations, Zhang\cite{Zh96} gave a proof for $d=3$ and Frantzikinakis and Kra\cite{FK05} for general $d$. Without those assumptions, the $L^{2}$-convergence of (\ref{1.1}) was established by Tao, as we have said above, it possesses four proofs. When $T_{1},T_{2},\dots,T_{d}$ belongs to nilpotent, it was proved by Miguel N.Walsh\cite{MW12}.

Although most people believe the existence of (\ref{1.1}) almost everywhere, the case in which one knows the answer are scarce. But we know that (\ref{1.1}) has a subsequence which converges almost everywhere. Meanwhile, the ergodic behaviour of diagonal measures is obtained, which is the main result of this paper. Additionally, on the torus $\mathbb{T}^{d}$ with special rotation, such as $R_{\alpha_{1},\cdots,\alpha_{d}}:\mathbb{T}^{d}\rightarrow \mathbb{T}^{d},$ with $1,\alpha_{1},\cdots,\alpha_{d}$  are rationally independent.  We can not only get the convergence of (\ref{1.1}) for every point in $\mathbb{T}^{d}$, but also get a beautiful result for ergodic behaviour. E. H. EL Abdalaoui considered the pointwise convergence of multiple ergodic averages in \cite{Ea17}. Here we give a different proof.

Before launching into the main result, we first reminder the reader some elements of measure theory and ergodic theorem in section 2. With sufficient preparation, in section 3, we give a proof of the ergodic behaviour of Theorem \ref{tao}, and give a classification of $T_{1},T_{2},\dots,T_{d}$, with case 1: all the $T_{i}$ are pairwise different, i.e. $T_{i}\neq T_{j}, i\neq j$. Case 2: there are $k$, with $1\leq k\leq d$, such that $T_{i_{1}}=T_{i_{2}}=\dots=T_{i_{k}}$. In section 4, we use the result in section 3 to the special case which the space is the torus $\mathbb{T}^{d}$, and transformations $R_{\alpha_{1}},\cdots,R_{\alpha_{d}}:\mathbb{T}\rightarrow \mathbb{T},$ satisfied that  $1,\alpha_{1},\cdots,\alpha_{d}$  are rationally independent. In section 5, we will give examples to show that each alternative in section 3 and section 4 really occurs.

\section{Preliminary}
\setcounter{equation}{0}
Let us first recall from \cite{Wa82}\cite{Ru04} some basic facts on measure theory and ergodic theorem.

\subsection{Measure Theory}
In this section, $X$ will be an arbitrary measure space with a positive measure $\mu$.
\bb{defn}\cite{Ru04}
Let $\mu$ be a positive measure on $X$. A sequence $\{f_{n}\}$ of measurable functions on $X$ is said to convergence in measure to the measurable function $f$ if for every $\epsilon>0$ there corresponds an $N$ such that
\be \lb{cm} \mu(\{x:|f_{n}(x)-f(x)|>\epsilon\})<\epsilon\ee
for all $n>N$.
\end{defn}
\bb{defn}\cite{Ru04}
If $1\leq p<\infty$ and if $f$ is a measurable function on $X$, define
\be  \|f\|_{p}=\{\int_{X}|f|^{p}\dmu\}^{\f 1 p}\ee
and let $L^{p}(\mu)$ consist of all $f$ for which
\be  \|f\|_{p}<\infty.\ee
We call $\|f\|_{p}$ the $L^{p}$-norm of $f$.
\end{defn}
If $f,f_{1},\cdots,f_{n},\cdots\in L^{p}(\mu)$ with $\lim_{n\to \infty}\|f_{n}-f\|_{p}=0$, we say that $\{f_{n}\}$ converges to $f$ in the mean of order $p$, or that $\{f_{n}\}$ is $L^{p}$-convergent to $f$.
\bb{thm}\lb{ae}
Assuming $\mu(X)<\infty$, we have the following statements:
\begin{enumerate}
  \item If $f_{n}\in L^{p}(\mu)$ and $\|f_{n}-f\|_{p}\rightarrow0$, then $f_{n}\rightarrow f$ in measure; here $1\leq p <\infty$.
  \item If $f_{n}\rightarrow f$ in measure, then $\{f_{n}\}$ has a subsequence $\{f_{n_{i}}(x)\}$ which converges to $f$ almost everywhere, i.e.
$$\lim_{i\to \infty}f_{n_{i}}(x)=f(x), a.e. x.$$
\end{enumerate}
\end{thm}

\subsection{Ergodic Theorem}
Let $(X,\mathcal{B},\mu)$ be a probability space and $T:X\rightarrow X$ be a measure-preserving transformation. $T$ is called ergodic if $T^{-1}B=B$ for $B\in\mathcal{B}$ satisfy $\mu(B)=0$ or $\mu(B)=1$. The next Theorem gives another form of the definition of ergodicity.

\bb{thm} \lb{er} {\rm \cite{Wa82}} Let $(X,\mathcal{B},\mu)$ be a probability space and let $T:X\rightarrow X$ be a measure-preserving transformation. Then $T$ is ergodic iff ~$\forall A,B \in \mathcal{B} $
\be \lb{mea} \frac{1}{n}\sum_{i=0}^{n-1} \mu(T^{-i}A\cap B)\rightarrow \mu(A)\mu(B). \ee
\end{thm}
To describe the proof of our main result, it will be convenient to reformulate Theorem \ref{er} in terms of functions.
\bb{cor}\lb{fun} Let $(X,\mathcal{B},\mu)$ be a probability space and let $T:X\rightarrow X$ be a measure-preserving transformation. Then $T$ is ergodic iff~ $\forall f,g:X\rightarrow \mathbb{R}, f,g\in L^{\infty}(X,\mathcal{B},\mu)$,
\be  \lim_{n\to \infty}\frac{1}{n}\sum_{i=0}^{n-1}\int f(T^{i}(x))g(x)\dmu=\int f(x)\dmu\int g(x)\dmu. \ee
\end{cor}
\Proof
Let us first prove the part $\Rightarrow$. Assume that $T$ is ergodic, from Theorem \ref{er}, we have
\be \lb{fun1} \lim_{n\to \infty}\frac{1}{n}\sum_{i=0}^{n-1}\int \chi_{A}(T^{i}(x))\chi_{B}(x)\dmu=\int \chi_{A}(x)\dmu\int \chi_{B}(x)\dmu, \ee
The functions in $L^{\infty}$ can be approximated, based on the function approximation theory, with some simple functions. we obtain the desired result (\ref{fun}).\\
Now let us prove the part $\Leftarrow$. Let $f=\chi_{A}, g=\chi_{B}$, we can get (\ref{mea}) easily. From Theorem \ref{er}, $T$ is ergodic.
\qed

The first major result in ergodic theory was proved in 1931 by G.D. Birkhoff. We will just state it.
\bb{thm}\lb{Birk}\cite{Wa82}
Suppose $T:(X,\mathcal{B},\mu)\rightarrow(X,\mathcal{B},\mu)$ is measure-preserving and $f\in L^{1}(\mu)$. Then $\disp \f 1 N \sum_{n=0}^{N-1} f(T^{n}x)$ converges almost everywhere to a function $f^{*}\in L^{1}(\mu)$. Also $f^{*}\circ T=f^{*}$ almost everywhere and if $\mu(X)<\infty$, then $\int f^{*}\dmu=\int f\dmu$.
\end{thm}

\bb{rem}\lb{ergodic}
If $T$ is ergodic then $f^{*}$ is constant almost everywhere and so if $\mu(X)<\infty$, $f^{*}=(1/ \mu(X))\int f\dmu$ almost everywhere. If $(X,\mathcal{B},\mu)$ is a probability space and $T$ is ergodic we have $\forall f\in L^{1}(\mu)$, $\lim_{N\to \infty}\frac{1}{N}\sum_{n=0}^{N-1}\int f(T^{n}(x))\dmu=\int f(x)\dmu$ almost everywhere.
\end{rem}

Unique ergodicity can get much more stronger behaviour in the ergodic theorem, now let me recall the definition of unique ergodicity and its strong behaviour of those averages in Theorem \ref{Birk}.
\bb{defn}\cite{Wa82}
A continuous transformation $T:X\rightarrow X$, where $X$ is a compact metrisable space, is called uniquely ergodic if there is only one $T$ invariant Borel probability measure on $X$.
\end{defn}

\bb{thm}\cite{Wa82}\lb{ue} Let $T:X\rightarrow X$ be a continuous transformation of a compact metrisable space $X$. The following are equivalent:\\
\begin{itemize}
\item For every $f\in C(X), \disp \f 1 N \sum_{n=0}^{N-1} f(T^{n}x)$ converges uniformly to a constant.
\item For every $f\in C(X), \disp \f 1 N \sum_{n=0}^{N-1} f(T^{n}x)$ converges pointwise to a constant.
\item There is a unique probability measure on $X$ which is invariant under $T$ such that for all $f\in C(X)$ and all $x\in X$,
    $$\f 1 N \sum_{n=0}^{N-1} f(T^{n}x)\rightarrow \int f\dmu.$$
\item $T$ is uniquely ergodic.
\end{itemize}

\end{thm}

When $X$ is a compact metric space, let $\B(X)$ be the Borel $\sigma$-algebra of $X$, $T:X\rightarrow X$ be a continuous transformation. We shall denote $\M(X),\M(X,T)$
$$\M(X)=\{\mu:\B(X)\rightarrow [0,1] \mid\mu(X)=1\},$$
$$\M(X,T)=\{\mu\in\M(X)\mid\mu(T^{-1}B)=\mu(B),B\in\B(X)\}.$$
 we know that
$\M(X)$ is convex and compact in the weak$^{*}$-topology\cite[Theorem 6.5]{Wa82}. $\M(X,T)$ is a convex and compact subset of $\M(X)$\cite[Theorem 6.5]{Wa82}.
\bb{cor} \lb{common measure} If $T_{1},T_{2},\cdots,T_{d}$ are pairwise commuting continuous maps of a compact space $X$ into itself, they possess a common invariant probability measure.
\end{cor}

\section{Ergodic Behaviour}\lb{ergodic behaviour}
Ergodic theory is the study of statistical properties of dynamical systems related to a measure on the underlying space of the dynamical system. The name comes from classical statistical mechanics, where the "ergodic hypothesis" asserts that, asymptotically, the time averages of an observable is equal to the space average. Now, we have known the "ergodic hypothesis" happens if the system is ergodic(see Remark \ref{ergodic}).

Before the main result of this paper is gave, let us introduce the concept of an irreducible dynamical system as a preliminary.
\bb{defn} The probability space $(X,\mathcal{B},T_{1},T_{2},\dots,T_{d},\mu)$, where $T_{1},T_{2},\dots,T_{d}:X\to X$ are
commuting invertible measure-preserving transformations, is called irreducible. If for any $f_{1},f_{2},\dots,f_{d}\in L^{\infty}(X,\mathcal{B},\mu)$,
$$A=\{x:\lim_{n\to \infty}\frac{1}{N} \sum_{n=0}^{N-1} f_{1}(T_{1}^{n}x)\cdots f_{d}(T_{d}^{n}x)<\alpha\},$$
$$B=\{x:\lim_{n\to \infty}\frac{1}{N} \sum_{n=0}^{N-1} f_{1}(T_{1}^{n+1}x)\cdots f_{d}(T_{d}^{n+1}x)<\alpha\},$$
$\forall\alpha \in \mathbb{R}$, we have $A=B$, and $A,B$ is the invariant set.
\end{defn}

\bb{thm} \lb{eh} Let $d\geq1$ be an integer. Assuming that $T_{1},T_{2},\dots,T_{d}:X\to X$ are
commuting invertible measure-preserving transformations of a irreducible space $(X,\mathcal{B},\mu)$, $T_{i}T_{j}^{-1}, i\neq j$
are ergodic. Then, for any $f_{1},f_{2},\dots,f_{d}\in L^{\infty}(X,\mathcal{B},\mu)$.
 $T_{i}\neq T_{j},i\neq j$,
\be  \lim_{k \to \oo}\disp \f 1 N_{k} \sum_{n=0}^{N_{k}-1} f_{1}(T_{1}^{n}x)\cdots f_{d}(T_{d}^{n}x)=
 \int_{X} f_{1}\dmu\cdots \int_{X} f_{d}\dmu \qq a.e..\ee
\end{thm}

\Proof
According to Theorem \ref{tao}, (\ref{1.1}) converges in $L^{2}$, According to Theorem \ref{ae}, (\ref{1.1}) has a subsequence which converges almost everywhere.\\
First of all, We shall prove that under the condition of irreducible, the subsequence converges to a constant almost everywhere. Assuming the contrary, we denote by $M$ the least upper bound of $f(x)$ over $X$ computed on neglecting a set of measure zero and analogously we denote by $m$ the greatest lower bound of the function $f(x)$ on neglecting a set of measure zero. From the assumption there follows $M>m$.\\
Let $\alpha$ satisfy the inequalities $m<\alpha<M$. We obtain
$$\mu\{p:f(p)<\alpha\}=\mu E_{\alpha}>0,$$
and
$$\mu(X-E_{\alpha})=\mu\{p:f(p)\geq\alpha\}>0,$$
we have a decomposition of $X$ into two sets of positive measure which contradicts the condition of irreducibility.\\
For $d=1$, the statement is the mean ergodic theorem. \\
For $d=2$, and almost every $x\in X$
\beaa  \lim_{k \to \oo}\disp \f 1 N_{k} \sum_{n=0}^{N_{k}-1} f_{1}(T_{1}^{n}x) f_{2}(T_{2}^{n}x)\EQ\int_{X}\lim_{k \to \oo}\disp \f 1 N_{k} \sum_{n=0}^{N_{k}-1} f_{1}(T_{1}^{n}x) f_{2}(T_{2}^{n}x)\dmu\\
\EQ \lim_{k \to \oo}\int_{X}\disp \f 1 N_{k} \sum_{n=0}^{N_{k}-1} f_{1}(T_{1}^{n}x)f_{2}(T_{2}^{n}x)\dmu\\
\EQ \lim_{k \to \oo}\int_{X}\disp \f 1 N_{k} \sum_{n=0}^{N_{k}-1} f_{1}(x)f_{2}((T_{2}T_{1}^{-1})^{n}x)\dmu\\
\EQ \lim_{k \to \oo}\disp \f 1 N_{k} \sum_{n=0}^{N_{k}-1}\int_{X} f_{1}(x)f_{2}((T_{2}T_{1}^{-1})^{n}x)\dmu\\
\EQ \int_{X} f_{1}\dmu\int_{X} f_{2}\dmu.\qq \mbox{(by Corollary \ref{fun})}
\eeaa
In the procedure of the proof, we used the Lebesgue's Dominated Convergence Theorem\cite[p. 26]{Ru04}.\\
Putting $f_{1}=\chi_{A},f_{2}=\chi_{B}$ in the situation $d=2$, then
$$
\lim_{k \to \oo}\disp \f 1 N_{k} \sum_{n=0}^{N_{k}-1} \chi_{A}(T_{1}^{n}x)\chi_{B}(T_{2}^{n}x)= \int_{X} \chi_{A}\dmu\int_{X} \chi_{B}\dmu \qq a.e..
$$
Multiplying by $\chi_{C}$,
$$
\lim_{k \to \oo}\disp \f 1 N_{k} \sum_{n=0}^{N_{k}-1} \chi_{T_{1}^{-n}A\cap T_{2}^{-n}B}\chi_{C}= \int_{\mathbb{T}} \chi_{A}\dmu\int_{\mathbb{T}} \chi_{B}\dmu\cdot\chi_{C}\qq a.e.,
$$
and the dominated convergence theorem implies
\be
\lim_{k \to \oo}\disp \f 1 N_{k} \sum_{n=0}^{N_{k}-1} \mu(T_{1}^{-n}A\cap T_{2}^{-n}B\cap C)= \mu(A)\mu(B)\mu(C).
\ee
For $d=3$, almost every $x\in X$
\begin{align*}
  &\lim_{k \to \oo}\disp \f 1 N_{k} \sum_{n=0}^{N_{k}-1} f_{1}(T_{1}^{n}x)f_{2}(T_{2}^{n}x)f_{3}(T_{3}^{n}x)\\
  =&\int_{X}\lim_{k \to \oo}\disp \f 1 N_{k} \sum_{n=0}^{N_{k}-1} f_{1}(T_{1}^{n}x)f_{2}(T_{2}^{n}x)f_{3}(T_{3}^{n}x)\dmu\\
 =& \lim_{k \to \oo}\int_{X}\disp \f 1 N_{k} \sum_{n=0}^{N_{k}-1} f_{1}(T_{1}^{n}x)f_{2}(T_{2}^{n}x)f_{3}(T_{3}^{n}x)\dmu\\
 =&\lim_{k \to \oo}\int_{X}\disp \f 1 N_{k} \sum_{n=0}^{N_{k}-1} f_{1}(x)f_{2}((T_{2}T_{1}^{-1})^{n}x)f_{3}((T_{3}T_{1}^{-1})^{n}x)\dmu\\
 =&\lim_{k \to \oo}\disp \f 1 N_{k} \sum_{n=0}^{N_{k}-1}\int_{X} f_{1}(x)f_{2}((T_{2}T_{1}^{-1})^{n}x)f_{3}((T_{3}T_{1}^{-1})^{n}x)\dmu\\
 =&\int_{X} f_{1}\dmu\int_{X} f_{2}\dmu\int_{X} f_{3}\dmu.
\end{align*}

In the same way, after $d$ steps, we obtain
$$ \lim_{k \to \oo}\disp \f 1 N_{k} \sum_{n=0}^{N_{k}-1} f_{1}(T_{1}^{n}x)\cdots f_{d}(T_{d}^{n}x)= \int_{X} f_{1}\dmu\cdots \int_{X} f_{d}\dmu.$$
\qed

\bb{cor} Let $d\geq1$ be an integer. Assuming that $T_{1},T_{2},\dots,T_{d}:X\to X$ are
commuting invertible measure-preserving transformations of a irreducible space $(X,\mathcal{B},\mu)$, with $T_{i_{1}}=T_{i_{2}}=\dots=T_{i_{k}}$, $1\leq k\leq d$. Without loss of generality, let $i_{1}=1,i_{2}=2,\dots,i_{k}=k$.  $T_{i}T_{j}^{-1}, i\neq j$
are ergodic. Then, for any $f_{1},f_{2},\dots,f_{d}\in L^{\infty}(X,\mathcal{B},\mu)$,
\be  \lim_{k \to \oo}\disp \f 1 N_{k} \sum_{n=0}^{N_{k}-1} f_{1}(T_{1}^{n}x)\cdots f_{d}(T_{d}^{n}x)=
 \int_{X}f_{1}f_{2}\cdots f_{k}\dmu\int_{X} f_{k+1}\dmu\cdots \int_{X} f_{d}\dmu \qq a.e..\ee
\end{cor}

\Proof
As the proof of Theorem\ref{eh}, the subsequence converges to a constant almost everywhere. \\
For $d=1$, the statement is the mean ergodic theorem. \\
For $d=2$, if $T_{1}\neq T_{2}$, it is conformity with Theorem\ref{eh}.\\
For $d=2$, and $T_{1}=T_{2}=T$, for almost every $x\in X$
\beaa  \lim_{k \to \oo}\disp \f 1 N_{k} \sum_{n=0}^{N_{k}-1} f_{1}(T_{1}^{n}x) f_{2}(T_{2}^{n}x)\EQ\int_{X}\lim_{k \to \oo}\disp \f 1 N_{k} \sum_{n=0}^{N_{k}-1} f_{1}(T^{n}x) f_{2}(T^{n}x)\dmu\\
\EQ\lim_{k \to \oo}\disp \f 1 N_{k} \sum_{n=0}^{N_{k}-1} \int_{X}f_{1}(T^{n}x) f_{2}(T^{n}x)\dmu\\
\EQ \int_{X} f_{1}f_{2}\dmu.\qq \mbox{($\mu$ is invariant)}
\eeaa
In the procedure of the proof, we used the Lebesgue's Dominated Convergence Theorem\cite[p. 26]{Ru04}.\\
For $d=3$, if $T_{1}\neq T_{2}\neq T_{3}$, it is conformity with Theorem\ref{eh}.\\
For $d=3$, and $T_{1}=T_{2}=T\neq T_{3}$, for almost every $x\in X$,
\beaa  \lim_{k \to \oo}\disp \f 1 N_{k} \sum_{n=0}^{N_{k}-1} f_{1}(T_{1}^{n}x)f_{2}(T_{2}^{n}x)f_{3}(T_{3}^{n}x)\EQ\int_{X}\lim_{k \to \oo}\disp \f 1 N_{k} \sum_{n=0}^{N_{k}-1} f_{1}(T^{n}x)f_{2}(T^{n}x)f_{3}(T_{3}^{n}x)\dmu\\
\EQ \lim_{k \to \oo}\disp \f 1 N_{k} \sum_{n=0}^{N_{k}-1}\int_{X} f_{1}(T^{n}x)f_{2}(T^{n}x)f_{3}(T_{3}^{n}x)\dmu\\
\EQ \lim_{k \to \oo}\disp \f 1 N_{k} \sum_{n=0}^{N_{k}-1}\int_{X} f_{1}(x)f_{2}(x)f_{3}((T_{3}T^{-1})^{n}x)\dmu\\
\EQ \int_{X} f_{1}f_{2}\dmu\int_{X} f_{3}\dmu.
\eeaa
For $d=3$, and $T_{1}=T_{2}=T_{3}=T$, for almost every $x\in X$,
\beaa  \lim_{k \to \oo}\disp \f 1 N_{k} \sum_{n=0}^{N_{k}-1} f_{1}(T_{1}^{n}x)f_{2}(T_{2}^{n}x)f_{3}(T_{3}^{n}x)\EQ\int_{X}\lim_{k \to \oo}\disp \f 1 N_{k} \sum_{n=0}^{N_{k}-1} f_{1}(T^{n}x)f_{2}(T^{n}x)f_{3}(T^{n}x)\dmu\\
\EQ \lim_{k \to \oo}\disp \f 1 N_{k} \sum_{n=0}^{N_{k}-1}\int_{X} f_{1}(T^{n}x)f_{2}(T^{n}x)f_{3}(T^{n}x)\dmu\\
\EQ \int_{X} f_{1}f_{2}f_{3}\dmu.
\eeaa
In the same way, after $d$ steps, we obtain
$$ \lim_{k \to \oo}\disp \f 1 N_{k} \sum_{n=0}^{N_{k}-1} f_{1}(T_{1}^{n}x)\cdots f_{d}(T_{d}^{n}x)=\int_{X}f_{1}f_{2}\cdots f_{k}\dmu\int_{X} f_{k+1}\dmu\cdots \int_{X} f_{d}\dmu \qq a.e..$$
\qed

\bb{rem}
If the more difficult question of pointwise almost everywhere convergence of (\ref{1.1}) is obtained, our methods can be used to the sequence as well.
\end{rem}

\section{Special case} \lb{special case}
\setcounter{equation}{0}

Consider on the torus $\mathbb{T}^{d}$ with special rotation, we can not only get the convergence of (\ref{1.1}) for every point in $\mathbb{T}^{d}$, although almost everywhere convergence of (\ref{1.1}) still don't known, but also get the result of Theorem \ref{eh}.

First of all, let me recall the definition of rationally independent.
\bb{defn}
The real numbers $1,\alpha_{1},\cdots,\alpha_{d}$ are rationally independent, if there is no $k_{0},k_{1},\cdots,k_{d}\in Z^{d+1}\backslash\{0\}$ such that $k_{0}+k_{1}\alpha_{1}+\cdots k_{d}\alpha_{d}=0$.
\end{defn}

The translation
$$R_{\alpha}=R_{\alpha_{1},\cdots,\alpha_{d}}(x_{1},\cdots,x_{d})=(x_{1}+\alpha_{1},\cdots,x_{d}+\alpha_{d})$$
induces the rotation $R_{\alpha}=R_{\alpha_{1},\cdots,\alpha_{d}}:\mathbb{T}^{d}\rightarrow \mathbb{T}^{d}.$

From next lemma, the rationally independent rotation $R_{\alpha}=R_{\alpha_{1},\cdots,\alpha_{d}}:\mathbb{T}^{d}\rightarrow \mathbb{T}^{d}$ is uniquely ergodic.
\bb{lem} \lb{unique-ergodicity} If $\alpha=(\alpha_{1},\cdots,\alpha_{d})$ with $1,\alpha_{1},\cdots,\alpha_{d}$ rationally independent, the Haar measure is the only probability measure which is invariant by $R_{\alpha}:\mathbb{T}^{d}\rightarrow \mathbb{T}^{d}.$
\end{lem}

Analogy Theorem \ref{eh}, the following statements satisfied:
\bb{thm} \lb{seh} Let $d\geq1$ be an integer. If $\alpha=(\alpha_{1},\cdots,\alpha_{d})$ with $1,\alpha_{1},\cdots,\alpha_{d}$ rationally independent, $R_{\alpha_{1}},\cdots,R_{\alpha_{d}}:\mathbb{T}\rightarrow \mathbb{T}$,  $f_{1},\cdots,f_{d}:\mathbb{T}\rightarrow \mathbb{R}$, $f_{1},f_{2},\dots,f_{d}\in C(\mathbb{T})$. Assuming $\mu$ is the common invariant probability measure of $R_{\alpha_{1}},\cdots,R_{\alpha_{d}}$.
Then
\be  \lim_{N \to \oo}\disp \f 1 N \sum_{n=0}^{N-1} f_{1}(R_{\alpha_{1}}^{n}x)\cdots f_{d}(R_{\alpha_{d}}^{n}x)=
 \int_{\mathbb{T}} f_{1}\dmu\cdots \int_{\mathbb{T}} f_{d}\dmu.
\ee
\end{thm}

\Proof  From Lemma \ref{unique-ergodicity}, for $\forall f \in C(\mathbb{T}^{d}),1\leq p <\infty$,
\be\lim_{N \to \oo}\disp \f 1 N \sum_{n=0}^{N-1} f(R_{\alpha}^{n}(y))=\int_{\mathbb{T}^{d}} f\dnu,\forall y\in\mathbb{T}^{d}, \ee
where $\nu=\mu\times\cdots\times\mu$.\\
Let $f(x_{1},x_{2},\cdots,x_{d})=f_{1}(x_{1})\cdot f_{2}(x_{2})\cdot\dots\cdot f_{d}(x_{d})$,
Then
\beaa  \lim_{N \to \oo}\disp \f 1 N \sum_{n=0}^{N-1} f_{1}(R_{\alpha_{1}}^{n}x)\cdots f_{d}(R_{\alpha_{d}}^{n}x)\EQ \lim_{N \to \oo}\disp \f 1 N \sum_{n=0}^{N-1} f(R_{\alpha}^{n}(x,\cdots,x))\\
\EQ \int_{\mathbb{T}^{d}} f\dnu\\
\EQ \int_{\mathbb{T}} f_{1}\dmu\cdots \int_{\mathbb{T}} f_{d}\dmu.
\eeaa
\qed

\bb{cor} Let $d\geq1$ be an integer. If $\alpha=(\alpha_{1},\cdots,\alpha_{d})$ with $1,\alpha_{1},\cdots,\alpha_{d}$ rationally independent, $R_{\alpha_{1}},\cdots,R_{\alpha_{d}},S:\mathbb{T}\rightarrow \mathbb{T}$, $f_{1},\cdots,f_{d},g:\mathbb{T}\rightarrow \mathbb{R}$, $f_{1},f_{2},\dots,f_{d}\in C(\mathbb{T})$, $g$ is a periodic function with $g(S^{k}x)=g(x)$. Assuming $\mu$ is the common invariant probability measure of $R_{\alpha_{1}},\cdots,R_{\alpha_{d}}$.
Then
\be  \lim_{N \to \oo}\disp \f 1 N \sum_{n=0}^{N-1} f_{1}(R_{\alpha_{1}}^{n}x)\cdots f_{d}(R_{\alpha_{d}}^{n}x)g(S^{n}x)=
 \int_{\mathbb{T}} f_{1}\dmu\cdots \int_{\mathbb{T}} f_{d}\dmu\f 1 k\sum_{r=0}^{k-1}g(S^{r}x).
\ee
\end{cor}
\Proof
When $N=kp$, then
\beaa  \lim_{p \to \oo}\disp \f 1 {kp} \sum_{n=0}^{kp-1} f_{1}(R_{\alpha_{1}}^{n}x)\cdots f_{d}(R_{\alpha_{d}}^{n}x)g(S^{n}x)\EQ \lim_{p \to \oo}\disp \f 1 {kp} \sum_{r=0}^{k-1}\sum_{i=0}^{p-1} f_{1}(R_{\alpha_{1}}^{ik+r}x)\cdots f_{d}(R_{\alpha_{d}}^{ik+r}x)g(S^{r}x)\\
\EQ\int_{\mathbb{T}} f_{1}\dmu\cdots \int_{\mathbb{T}} f_{d}\dmu \f 1 k\sum_{r=0}^{k-1}g(S^{r}x).
\eeaa
When $N=kp+m$, $1\leq m<k$, then
\begin{equation}
\begin{aligned}
&\lim_{p \to \oo}\disp \f 1 {kp+m} \sum_{n=0}^{kp+m-1} f_{1}(R_{\alpha_{1}}^{n}x)\cdots f_{d}(R_{\alpha_{d}}^{n}x)g(S^{n}x) \\
&=\lim_{p \to \oo}\disp \f 1 {kp+m} \{\sum_{r=0}^{k-1}\sum_{i=0}^{p-1} f_{1}(R_{\alpha_{1}}^{ik+r}x)\cdots f_{d}(R_{\alpha_{d}}^{ik+r}x)g(S^{r}x)\\
&~~+\sum_{j=1}^{m}f_{1}(R_{\alpha_{1}}^{pk+j}x)\cdots f_{d}(R_{\alpha_{d}}^{pk+j}x)g(S^{j}x)\}\\
&=\lim_{p \to \oo}\disp \f {kp} {kp+m}\f 1 {kp} \sum_{r=0}^{k-1}\sum_{i=0}^{p-1} f_{1}(R_{\alpha_{1}}^{ik+r}x)\cdots f_{d}(R_{\alpha_{d}}^{ik+r}x)g(S^{r}x)\\
&~~+\lim_{p \to \oo}\disp \f 1 {kp} \sum_{j=1}^{m}f_{1}(R_{\alpha_{1}}^{pk+j}x)\cdots f_{d}(R_{\alpha_{d}}^{pk+j}x)g(S^{j}x)\\
&=\int_{\mathbb{T}} f_{1}\dmu\cdots \int_{\mathbb{T}} f_{d}\dmu \f 1 k\sum_{r=0}^{k-1}g(S^{r}x).
\end{aligned}
\end{equation}
\qed

\section{Examples}
\setcounter{equation}{0}
Now we give examples to show that each alternative in section \ref{ergodic behaviour} and \ref{special case} really occurs. Although in the uniquely ergodic theorem, we need the function $f$ is continuous. In our examples, we can get the result with function is not continuous.

\bb{prop}
Let $\{x\}$ denote the decimal part of $x$, $f(x):S^{1}\rightarrow \mathbb{R}$, $f(x)=\{x\}$, giving measure-preserving transformation $T:S^{1}\rightarrow S^{1}$ is ergodic, $\mu$ is the Haar measure. Then
\be  \lim_{N \to \oo}\disp \f 1 N \sum_{n=0}^{N-1} f(T^{n}x)=\f 1 2.\ee
\end{prop}

\Proof
 Let $g(x),h(x)$ be
\begin{equation*}
g(x)=\left\{
\begin{array}{cc}
                    f(x) & 0\leq x\leq1-\f 1 m,\\
                    (1-m)(x-1) & 1-\f 1 m\leq x\leq1.
\end{array}\right.
\end{equation*}

\begin{equation*}
h(x)=\left\{
\begin{array}{cc}
                    (1-m)x+1 & 0\leq x\leq \f 1 m,\\
                    f(x) & \f 1 m\leq x\leq1.
\end{array}\right.
\end{equation*}
Obviously, $g(x),h(x)$ have the following properties:
\begin{enumerate}
  \item $g(x),h(x)$ are continuous functions,
  \item\lb{bd} $g(x)\leq f(x)\leq h(x),$\quad $\forall x\in S^{1},$
  \item $\lim_{m \to \oo}g(x)=\lim_{m \to \oo}h(x)=f(x),$
  \item $\int_{S^{1}} g(x)\dmu=(m-1)/(2m),$
  \item $\int_{S^{1}} h(x)\dmu=(m+1)/(2m).$
\end{enumerate}
By Theorem \ref{ue},
\be  \lim_{N \to \oo}\disp \f 1 N \sum_{n=0}^{N-1} g(T^{n}x)=\int_{S^{1}} g(x)\dmu.\ee
\be  \lim_{N \to \oo}\disp \f 1 N \sum_{n=0}^{N-1} h(T^{n}x)=\int_{S^{1}} h(x)\dmu.\ee
From property \ref{bd}, we have
$$\lim_{N \to \oo}\disp \f 1 N \sum_{n=0}^{N-1} g(T^{n}x)\leq\liminf_{N \to \oo}\disp \f 1 N \sum_{n=0}^{N-1} f(T^{n}x)\leq\limsup_{N \to \oo}\disp \f 1 N \sum_{n=0}^{N-1} f(T^{n}x)\leq\lim_{N \to \oo}\disp \f 1 N \sum_{n=0}^{N-1} h(T^{n}x).$$
According to all the properties and equations above, we have
$$\lim_{N \to \oo}\disp \f 1 N \sum_{n=0}^{N-1} f(T^{n}x)=\f 1 2.$$
\qed

\bb{ex}
Let $\{x\}$ denote the decimal part of $x$.
\begin{enumerate}
         \item \lb{ad1} $\lim_{N \to \oo}\disp \f 1 N \sum_{n=0}^{N-1} \{x+\sqrt{2}n\}\{x+\sqrt{3}n\}=\f 1 4.$
         \item \lb{ad2} $\lim_{N \to \oo}\disp \f 1 N \sum_{n=0}^{N-1} \{x+\sqrt{2}n\}\{x+\sqrt{2}n\}=\f 1 3.$
         \item \lb{ad3} $\lim_{N \to \oo}\disp \f 1 N \sum_{n=0}^{N-1} \{x+\sqrt{2}n\}\{x+\f n k\}=\f 1 {2k} \{kx\}+\f {k-1} {4k}.$
\end{enumerate}
\end{ex}

\Proof Let $f_{1}, f_{2}:S^{1}\rightarrow \mathbb{R}$, with $f_{1}(x)=f_{2}(x)=\{x\}$
\begin{align*}
  T_{\sqrt{2}}:S^{1}&\rightarrow S^{1}\\
              x    &\mapsto\sqrt{2}+x~mod~1.\\
  T_{\sqrt{3}}:S^{1}&\rightarrow S^{1}\\
              x    &\mapsto\sqrt{3}+x~mod~1.
\end{align*}

\bu Proof of Example \ref{ad1}
\beaa  \lim_{N \to \oo}\disp \f 1 N \sum_{n=0}^{N-1} \{x+\sqrt{2}n\}\{x+\sqrt{3}n\}\EQ \lim_{N \to \oo}\disp \f 1 N \sum_{n=0}^{N-1} f_{1}(T_{\sqrt{2}}^{n}x)f_{2}(T_{\sqrt{3}}^{n}x)\\
\EQ \int_{S^{1}} f_{1}\dmu\int_{S^{1}} f_{2}\dmu\\
\EQ \f 1 2\times \f 1 2=\f 1 4.
\eeaa

\bu Proof of Example \ref{ad2}
\beaa  \lim_{N \to \oo}\disp \f 1 N \sum_{n=0}^{N-1} \{x+\sqrt{2n}\}\{x+\sqrt{2n}\}\EQ \lim_{N \to \oo}\disp \f 1 N \sum_{n=0}^{N-1} f_{1}(T_{\sqrt{2}}^{n}x)f_{2}(T_{\sqrt{2}}^{n}x)\\
\EQ \int_{S^{1}} f_{1}^{2}\dmu\\
\EQ \int_{S^{1}} x^{2}\dmu\\
\EQ \f 1 3.
\eeaa

\bu Proof of Example \ref{ad3}

\noindent{\bf Assertion}:\be \lb{final}
\sum_{r=0}^{k-1}\{x+\f r k\}=\{kx\}+\f {k-1} 2,~~ x\in[0,1].
\ee

Let me prove the assertion firstly:

\Proof
 When $\f i k\leq x\leq\f {i+1} k,~~i=0,1,\cdots,k-1$. We have $x+\f r k\leq1,~r=0,1,\cdots,k-i-1,$ and $1\leq x+\f r k\leq2,~r=k-i,k-i+1,\cdots,k-1,$ Then
  \begin{align*}
    \sum_{r=0}^{k-1}\{x+\f r k\}=&\sum_{r=0}^{k-i-1}(x+\f r k)+\sum_{r=k-i}^{k-1}(x+\f r k-1) \\
                                =&kx+\f {k-1} 2-i\\
                                =&\{kx\}+\f {k-1} 2.
    \end{align*}

\qed

Let $g:S^{1}\rightarrow \mathbb{R}$, with $g(x)=\{x\}$,
\begin{align*}
  T:S^{1}&\rightarrow S^{1}\\
    x    &\mapsto x+\f 1 k ~mod~1.
\end{align*}
\beaa  \lim_{N \to \oo}\disp \f 1 N \sum_{n=0}^{N-1} \{x+\sqrt{2}n\}\{x+\f n k\}\EQ \lim_{N \to \oo}\disp \f 1 N \sum_{n=0}^{N-1}f_{1}(T_{\sqrt{2}}^{n}x)g(T^{n}x)\\
\EQ \int_{S^{1}} f_{1}\dmu\cdot\f 1 k\sum_{r=0}^{k-1}g(T^{r}x)\\
\EQ\f 1 {2k}\sum_{r=0}^{k-1}\{x+\f r k\}\\
\EQ \f 1 {2k} \{kx\}+\f {k-1} {4k}.
\eeaa

\noindent{\bf Acknowledgements}:The authors thank to Xu Xia for valuable suggestions.

\end{document}